\documentclass[a4paper,fleqn]{cas-sc}

\newtheorem{theorem}{Theorem}[section]

\newtheorem{definition}{Definition}[section]

\newtheorem{lemma}{Lemma}
\newtheorem{assumption}{Assumption}
\newtheorem{algorithm}{algorithm}
\newtheorem{remark}{Remark}

\newproof{pf}{Proof}

\usepackage{graphics} 
\usepackage{amsmath} 
\usepackage{amssymb}  
\usepackage{enumitem}   
\usepackage{hyperref}
\usepackage{algorithm,algorithmic}
\usepackage[authoryear,longnamesfirst]{natbib}
\usepackage{caption}
\usepackage{subcaption}

\title{\LARGE \bf
Data-driven distributionally robust MPC using the Wasserstein metric
}

\begin{document}

\let\WriteBookmarks\relax
\def\floatpagepagefraction{1}
\def\textpagefraction{.001}
\shorttitle{Data-driven distributionally robust MPC using the Wasserstein metric}
\shortauthors{Zhengang Zhong, Ehecatl Antonio del Rio-Chanona and Panagiotis Petsagkourakis}

\author%
{Zhengang Zhong}
\fnmark[1]

\author%
{Ehecatl Antonio del Rio-Chanona}\cormark[1]
\fnmark[1]

\author%
{Panagiotis Petsagkourakis}
\fnmark[2]\cormark[1]

\cortext[cor1]{Corresponding author.}

\fntext[fn1]{Zhengang Zhong and Ehecatl Antonio del Rio-Chanona are with Centre for Process Systems Engineering, Imperial College London, London, United Kingdom 
        {\tt\small  z.zhong20@imperial.ac.uk, a.del-rio-chanona@imperial.ac.uk}}
\fntext[fn2]{Panagiotis Petsagkourakis is with Centre for Process Systems Engineering, University College London, London, United Kingdom
        {\tt\small p.petsagkourakis@ucl.ac.uk}}%
        
\nonumnote{This project has also received funding from the EPSRC project EP/R032807/1.
  }

\begin{abstract}
A data-driven MPC scheme is proposed to safely control constrained stochastic linear systems using distributionally robust optimization.
Distributionally robust constraints based on the Wasserstein metric are imposed to bound the state constraint violations in the presence of process disturbance.
A feedback control law is solved to guarantee that the predicted states comply with constraints. The stochastic constraints are satisfied with regard to the worst-case distribution within the Wasserstein ball centered at their discrete empirical probability distribution.
The resulting distributionally robust MPC framework is computationally tractable and efficient, as well as recursively feasible. The innovation of this approach is that all the information about the uncertainty can be determined empirically from the data. The effectiveness of the proposed scheme is demonstrated through numerical case studies.
\end{abstract}


\maketitle

\section{INTRODUCTION}

Model predictive control (MPC) has demonstrated remarkable success due to its ability to handle multivariate dynamics and constraints\cite{dobos2009dynamic, TOURETZKY20141292}.
MPC solves an open-loop optimal control problem at each sampling time based on a nominal model to decide a finite sequence of control actions from which the first element of the sequence is implemented \cite{MPC_book}.

In the context of control under uncertainty, two important methodologies arise to guarantee constraint satisfaction: Robust MPC (RMPC) \cite{RMPC_book_Houska} and stochastic MPC (SMPC) \cite{arcari2020dual, HEWING2020109095}.
The former addresses the receding horizon optimal control problem for uncertain systems in a deterministic fashion by assuming bounded uncertainties and providing solutions for the worst case-scenario.
Some important approaches to RMPC are min-max optimization \cite{campo1987robust} and tube-based MPC \cite{mayne2011tube}.
However, some worst-case scenarios are unlikely to work in practice; as resulting control designs tend to be over conservative or even infeasible \cite{garatti2013modulating}. In addition to conservativeness, if the assumed uncertainty sets are inaccurate, the controller may have poor performance.
To address this, approaches have been proposed  to reduce conservativeness in the context of RMPC, e.g. \cite{LUCIA201269, LUCIA20141904, 9145693}.

In contrast to RMPC, SMPC solves a stochastic optimization problem by assuming distributional information of the uncertainty \cite{Schwarm1999}.
Chance constraints on SMPC reduce the inherent conservativeness of robust MPC via the trade-off between constraint satisfaction and closed-loop performance \cite{mesbah2014stochastic}.
However, deviation of the assumed distribution from the true one caused by poor assumptions or limited available data may result in sub-optimality, infeasibility and unwanted behavior of the system \cite{nilim2005robust}. 

To overcome the conservativeness of RMPC and the distributional mismatch of SMPC, we explore a distributionally robust optimization approach - distributionally robust MPC.
In distributionally robust optimization (DRO), a variant of the stochastic optimization is explored where the probability distribution is itself uncertain.
DRO minimizes an expected loss function, where the expectation comes from the worst-case probability distribution of an ambiguity set.

Just as wrong assumptions on the distribution of the uncertainty can be detrimental to the objective's performance on an MPC scheme, chance constraints can also be affected by this mismatch, incurring in severe violations.
As a counterpart to chance constraints, distributionally robust chance constraints assume the actual distribution of uncertain variables belongs to an ambiguity set. This ambiguity set contains all distributions with a predefined characteristic (e.g. first or second moments), and such an ambiguity set can be computed from historical data.
Distributionally robust constraints have direct connection to the constraints incorporated in the classical paradigms of RMPC and SMPC \cite{rahimian2019distributionally}. Given the aforementioned benefits, distributionally robust constraints \cite{wiesemann2014distributionally} are considered in this work. We characterize the ambiguity set as a Wasserstein ball around an empirical distribution with a radius defined by the Wasserstein metric \cite{chen2018data}. 

The Wasserstein ambiguity set has received increasing attention in distributionally robust optimization due to its modeling flexibility, finite-sample guarantee, and tractable reformulation into convex programs \cite{esfahani2018data}.
In contrast to a Wasserstein ambiguity set, other ambiguity sets do not enjoy these properties. For example, approaches leveraging moment-based ambiguity do not have finite-sample guarantees \cite{luo2019decomposition} and ambiguity sets using $\phi$-divergence such as Kullback–Leibler divergence typically contain irrelevant distributions \cite{gao2016distributionally}.
These drawbacks motivate our use of a Wasserstein ambiguity set. 

Distributionally robust optimization in control problems has been studied with regard to different formulations of objective functions \cite{You2020}.
In the setting of multi-stage stochastic optimal power flow problem, a framework is proposed to solve multi-stage feedback control policies with respect a
Conditional Value at Risk (CVaR) objective \cite{guo2018data}.
A control policy for wind power ramp management is solved via dynamic programming by reformulating the distributionally robust value function with a tractable form: a convex piecewise linear ramp penalty function  \cite{yang2019data}.
For linear quadratic problems, a deterministic stationary policy is determined by solving an data-irrelevant discrete algebraic Riccati equation \cite{9222209}.
The control policies for both finite horizon optimal control problem with expected quadratic objective and infinite horizon optimal control problem with an average cost can be determined concerning distributionally robust CVaR constraints modeled by moment-based ambiguity sets \cite{van2015distributionally}.
Recently, a data-driven distributionally robust MPC with a moment-based ambiguity set for quadratic objective function under multi-stage risk measures was proposed in \cite{coppens2021data}.

\subsection{Main Contribution}
Compared to existing related studies which consider linear objective function or moment-based ambiguity sets, our work uniquely considers a finite-horizon control problem with distributionally robust constraints constructed by Wasserstein ambiguity set.
The main contributions of this paper are summarized as follows:
1. A model predictive controller with distributionally robust quadratic objective and distributionally robust chance constraints is proposed to determine purified-output-based (POB) affine control laws.
2.  A practical Algorithm is proposed, which results in a tractable conic optimization problem, and all the information about the uncertainty can be determined empirically from the data. 
3.  A proof for recursive feasibility is provided, and finite sample guarantee of chance constraints is illustrated via a case study.

\subsection{Notation}
Let $x_{[k:k+N]}$ denotes the concatenated state vector $[x_k^{\top}, x_{k+1}^{\top},\dots,x_N^{\top}]^{\top}$ and $[x_{[k: k+N]}]_{i}$ be the i-th entry of the vector. 
$\mathbb{D}(x_k)$ denotes the set of feasible state trajectory within the prediction horizon $N$ with the initial state $x_k$, and $\mathbb{I}_m$ is an index set indicating the $m$-th state is constrained $\forall m \in \mathbb{I}_m$. Let $\operatorname{Tr}( \cdot )$ be the trace operator.
We denote by $\mathbb{S}_{+}^{n}$ and $\mathbb{S}_{++}^{n}$ the sets of all positive semidefinite and positive definite symmetric matrices in $\mathbb{R}^{n \times n}$, respectively.
The diagonal concatenation of two matrices $A$ and $B$ is denoted by $\operatorname{diag}(A, B)$. $A_{i,j}$ is the entry of i-th row and j-th column in a 2D matrix and $[A]_{i,j:k}$ is the row vector of i-th row and j-th to k-th columns in matrix $A$. $e_j$ denote a column vector with all entries $0$ expect $j$-th entry equal to $1$ and $\mathbb{e}_j$ selects state vector at the sampling time $j$ from the stacked state vector.

All random vectors are defined as measurable functions on an abstract probability space $\left(\Omega, X, \mathbb{P}\right),$ where $\Omega$ is referred to as the sample space, $X$ represents the $\sigma$-algebra of events, and $\mathbb{P}$ denotes the true but unknown probability measure.
We denote by $\delta_\xi$ the Dirac distribution concentrating unit mass at $\xi$ and by $\delta_k$ the state difference of nominal and disturbed system at sampling time $k$. The $N$-fold product of a distribution $\mathbb{P}$ on $\Xi$ is denoted by $\mathbb{P}^{N}$, which represents a distribution on the Cartesian product space $\Xi^{N}$. $\mathcal{M}(\Xi)$ is the space of all probability distributions $\mathbb{Q}$ supported on $\Xi$ with finite expected norm. The training data set comprising $N_s$ samples is denoted by $\widehat{\Xi}_{N_s}:=\left\{\widehat{\xi}_{i}\right\}_{i<N_s} \subseteq \Xi$.
%
%
\subsection{Organization}
The remainder of this paper is organized as follows.
The problem formulation to determine the POB affine control laws for disturbed systems is introduced in Section II.
Preliminaries on distributionally robust control embracing the Wasserstein metric and the corresponding optimization problem are covered in Section III.
The main results including the tractable reformulation of the distributionally robust optimization problem, and finite sample guarantee are discussed in Section IV along with a practical Algorithm.
Simulation experiments for case studies mass spring system and inverted pendulum are illustrated in Section V, and the results are also analyzed.
Conclusions are summarized in Section VI.
%
\section{PROBLEM FORMULATION}
\label{section:problem_fomulation}
In this section we derive POB affine control laws \cite{ben2009robust} for a discrete-time linear time-invariant (LTI) dynamical system with additive disturbance. We consider a discrete-time LTI system at time $k$ 
\begin{equation}
    \begin{aligned}
x_{0} &=x \\
x_{k+1} &=A x_{k}+B u_{k}+C w_{k} \\
y_{k} &=D x_{k}+E w_{k},
\end{aligned} 
 \label{eq:DLTI} \tag{1}
\end{equation}
where state $x_{k} \in \mathbb{R}^{n_{x}}$, input $u_{k} \in \mathbb{R}^{n_{u}}$, output $y_{k} \in \mathbb{R}^{n_{y}}$ and the disturbance $w_{k} \in \mathbb{R}^{n_{w}}$.
Both process noise and measurement noise are modelled via matrices $C$ and $E$.
Our design target is to enable the closed-loop system of \eqref{eq:DLTI} to meet prescribed requirements. One of the requirements is to satisfy polyhedral state constraints $C_p x_{[k:k+N]}  \leq D_p$ within the prediction horizon $N$. 

The POB affine control laws are derived based on the discrepancy between the disturbed system and its corresponding nominal system (see Definition \ref{POB}).
Without loss of generality, we assume the equilibrium point is at the origin.
\begin{definition}[Model]
Given a disturbed system in the form of \eqref{eq:DLTI}, we define the corresponding nominal system initialized at the equilibrium point and not disturbed by exogenous inputs as the \emph{model}
\begin{equation}
\begin{aligned}
\widehat{x}_{0} &=0 \\
\widehat{x}_{k+1} &=A \widehat{x}_{k}+B u_{k} \\
\hat{y}_{k} &=D \widehat{x}_{k}.
\end{aligned}
\label{eq:model}
\tag{2}
\end{equation}
\end{definition}

The open-loop state difference between the model and the disturbed system $\delta_{k}=x_{k}-\widehat{x}_{k}$ evolves according to the system matrices and disturbances.
The accumulated influence of disturbances can be measured via purified outputs $v_{t}=D \delta_{t}+E w_{t} = y_{t}-\widehat{y}_{t}$ where $t \in \left[k, k+N-1\right]$. This allows us to consider POB affine control laws based on the disturbance and inputs history \cite{Ben-Tal2006}.

Hence, we solve the following problem:
Given a system's matrices, a horizon $N$, initial state $x_0 = x$ and sampled disturbance data points, we determine the control law that, by leveraging a tractable data-driven optimization problem, can steer the system to a desired equilibrium state while satisfying prescribed chance constraints $\mathbb{E}\left[C_{p} x_{[k: k+N]} \leq D_{p} \right] \ge 1-\beta$. The parameterized affine control laws will be defined as following.

\begin{definition}[POB Affine Control Laws]\label{POB}
At sampling time $t$, given purified outputs from $k$ to $t$, we define the \emph{POB affine control laws} as 
\begin{equation}
u_{t}=h_{t}+\sum_{\tau=k}^{t} H_{t, \tau} v_{\tau} 
\label{eq:POB_affine} \tag{3}  
\end{equation}
with $t \in \left[k, k+N - 1\right]$. 
\end{definition}
Note that the definition of the POB affine control law above is equivalent to an affine control law that is only dependent on initial state and disturbance sequence or only dependent on the outputs of the disturbed system.
\begin{lemma}[Equivalent Control Laws]
\label{thm:1}
For every POB affine control law in the form of \eqref{eq:POB_affine}, there exists a control law resulting in exactly the same closed-loop state-control trajectories dependent only on:
\begin{enumerate}[label=(\roman*)]
    \item initial state and disturbance sequences \label{thm:11}
    \item outputs of the disturbed system. \label{thm:12}
\end{enumerate}
\end{lemma}
\begin{pf} \ref{thm:11} follows directly from the definition of purified output. We reformulate \eqref{eq:POB_affine} by inserting the accumulated state discrepancy and the current disturbance. Then we acquire 
\begin{equation}
v_{t} = \sum_{\tau = 0}^{t} DA^{t-\tau}C^{\tau}x_0 + Ew_t 
\tag{4} \label{eq:affine_init_dist}
\end{equation}
Now we can show that \ref{thm:11} holds after inserting \eqref{eq:affine_init_dist} into \eqref{eq:POB_affine}.
The proof of \ref{thm:12} follows a similar reasoning as that in \cite[Theorem 14.4.1.]{ben2009robust}.
\end{pf}


In an effort to simplify notation when constructing the optimization problem to determine the control law, we derive a compact form of the dynamical system. The dynamics of the linear system over the finite horizon $N$ can be written as $x_{[k:k+N]}=A_{x} x + B_{x} u_{[k:k+N-1]} + C_{x} w_{[k:k+N-1]}$ and the corresponding measurements as $y_{[k:k+N-1]}=A_{y} x + B_{y} u_{[k:k+N-1]} + \left(C_{y} + E_{y}\right) w_{[k:k+N-1]}$, where $\small
A_{y}=\left[\begin{array}{l}
DA^{0} \\ \vdots \\DA^{N-1}
\end{array}\right]
$,
$\tiny
B_{y}=\left[\begin{array}{ccccc}
0 & 0  & \ldots & 0 & 0\\
DA^{0} B & 0  & \ldots & 0 & 0\\
\vdots & \ddots & \ddots  & \vdots & \vdots\\
\vdots & \vdots & \ddots  & 0 & 0\\
DA^{N-2} B & DA^{N-3} B  & \ldots & DA^{0} B & 0
\end{array}\right]
$ and $\small
E_{y} = \left[
  \begin{array}{ccc}
    E & &  \\
    & \ddots &  \\
    & &  E
  \end{array} 
  \right]
$. $A_x$, $B_x$ and $C_x$ are derived in a same way as $A_y$ and $B_y$ by formulating the stacked state in initial state, stacked input and stacked disturbance.
Consider now the inputs characterized by the POB affine control law, we derive $u_{[k:k+N-1]} = H_N(\tilde{C}_y + \tilde{E}_y) \tilde{w}_{[k:k+N-1]} = \tilde{H}_N\tilde{w}_{[k:k+N-1]} $, where $\tilde{w}_{[k:k+N-1]} =  \small \begin{bmatrix}1 & w_{[k:k+N-1]}^{\top}\end{bmatrix}^{\top} $ is the extended disturbance vector, $\tilde{C}_y = \begin{bmatrix} 0 & 0 \\ A_y x_0 & C_y \end{bmatrix}$ and $\tilde{C}_E = \begin{bmatrix} 1 & 0 \\ 0 & E_y \end{bmatrix}$. Furthermore, we use $H_N$ with subscript $N$ to denote the control law for horizon $N$. Finally, this allows us to write the stacked state vector as a linear matrix equation 
\begin{equation}
x_{[k:k+N]}=\left(B_{x} \tilde{H}_N + \tilde{C}_{x}\right) \tilde{w}_{[k:k+N-1]} \tag{5} \label{eq:stacked_state},
\end{equation}
where $\tilde{C}_x = \begin{bmatrix} A_x x_0 & C_x \end{bmatrix}$. Our goal is then to determine the control law $\tilde{H}_N$ which steers the system to the origin by minimizing the distributionally robust quadratic objective whilst guaranteeing distributionally robust constraint satisfaction.

\section{DISTRIBUTIONALLY ROBUST MPC}

\subsection{Ambiguity Sets and the Wasserstein metric}
Distributionally robust optimization is an optimization model where limited information about the underlying probability distribution of the random parameters in a stochastic model is assumed.
To model distributional uncertainty, we characterize the partial information about the true distribution $\mathbb{P}$ by a set of probability measures on $(\Omega, X)$. This set is termed as the \emph{ambiguity set} \cite{wiesemann2014distributionally}.
In this paper, we focus on an ambiguity set specified by a discrepancy model \cite{rahimian2019distributionally} wherein the distance on the probability distribution space is characterized by the Wasserstein metric.
The Wasserstein metric defines the distance between all probability distributions $\mathbb{Q}$ supported on $\Xi$ with finite $p$-moment $\int_{\Xi}\|\xi\|^{p} \mathbb{Q}(d \xi)<\infty$. 
\begin{definition}[Wasserstein Metric \cite{piccoli2014generalized}]
\! The \emph{Wasserstein metric} of order $p \ge 1$ is defined as $d_w: \mathcal{M}(\Xi) \times \mathcal{M}(\Xi) \rightarrow \mathbb{R}$ for all distribution $\mathbb{Q}_1, \mathbb{Q}_2 \in \mathcal{M}(\Xi)$ and arbitrary norm on $\mathbb{R}^{n_\xi}$:
\begin{equation}
d_{\mathrm{W}}\left(\mathbb{Q}_{1}, \mathbb{Q}_{2}\right):=\inf_{\Pi} \left\{\left(\int_{\Xi^{2}}\left\|\xi_{1}-\xi_{2}\right\|^{p} \Pi\left(\mathrm{d} \xi_{1}, \mathrm{~d} \xi_{2}\right)\right)^{1/p}\right.
\tag{6}    \label{eq:wasserstein_metric}
\end{equation}
where $\Pi$ is a joint distribution of $\xi_{1}$ and $\xi_{2}$ with marginals $\mathbb{Q}_{1}$ and $\mathbb{Q}_{2}$ respectively.
\end{definition}

Specifically, we define an ambiguity set centered at the empirical distribution leveraging the Wasserstein metric
\begin{equation}
\mathbb{B}_{\varepsilon}\left(\widehat{\mathbb{P}}_{N_s}\right):=\left\{\mathbb{Q} \in \mathcal{M}(\Xi): d_{\mathrm{W}}\left(\widehat{\mathbb{P}}_{N_s}, \mathbb{Q}\right) \leq \varepsilon\right\}
\tag{7} \label{eq:Wasserstein_ball}
\end{equation}
which specifies the Wasserstein ball with radius $\varepsilon>0$ around the discrete empirical probability distribution $\widehat{\mathbb{P}}_{N_s}$. The empirical probability distribution  $\widehat{\mathrm{P}}_{N_s}:=\frac{1}{N_s} \sum_{i=1}^{N_s} \delta_{\widehat{\xi}_{i}}$ is the uniform distribution on the training data set $\widehat{\Xi}_{N_s}:=\left\{\widehat{\xi}_{i}\right\}_{i \leq N_s} \subseteq \Xi$. $\delta_{\widehat{\xi}_{i}}$ is the Dirac distribution which concentrates unit mass at $\widehat{\xi}_{i} \in \mathbb{R}^{n_{\xi}}$.
The radius $\varepsilon$ should be large enough such that the ball contains the true distribution $\mathbb{P}$ with high fidelity, but not unnecessarily large to hedge against over-conservative solutions. The impact of the ball radius is illustrated and discussed in Section V.

\subsection{Data-Based Distributionally Robust MPC}
We consider the optimal control problem for the system \eqref{eq:stacked_state} enforcing distributionally robust constraints to be satisfied, i.e. 
\begin{equation}
\sup _{\mathbb{Q} \in \mathbb{B}_{\varepsilon}\left(\widehat{\mathbb{P}}_{N_s}\right)} \mathbb{E}^{\mathbb{Q}}[\ell(\xi,H_N)] \le U, 
\tag{8} \label{eq:DR_upperbound}
\end{equation}
where $\ell$ is a function representing state constraints in the polyhedral feasible space, dependent on affine control laws and disturbances, and $U$ is a prescribed bound. Distributionally robust constraints in a stochastic setting can take the information about the probability distribution into account via \eqref{eq:Wasserstein_ball} such that the prescribed state constraints hold in expectation for the worst-case distribution within the ball \eqref{eq:Wasserstein_ball}.

Our aim is to find an admissible affine control law with respect to the distributionally robust constraints whilst minimizing a distributionally robust objective $J_N$. We characterize the objective function as a discounted sum of quadratic stage costs
\begin{equation}
J_{N}\left(x, H_N\right):=\inf_{\mathbb{Q} \in \mathbb{B}_{\varepsilon}\left(\widehat{\mathbb{P}}_{N_s}\right)} \mathbb{E}_{\mathbb{Q}}\left\{\sum_{t=k}^{k+N-1} \beta^{t}\left[\boldsymbol{x}_{t}^{\top} Q \boldsymbol{x}_{t}+\boldsymbol{u}_{t}^{\top} R \boldsymbol{u}_{t}\right]  + \beta^{N} \boldsymbol{x}_{k+N}^{\top} Q_{f} \boldsymbol{x}_{k+N}\right\},
\tag{9} \label{eq:disc_quadratic_cost}
\end{equation}
with $\beta \in (0,1]$ as discount factor. It is further assumed that $Q, Q_{f} \in \mathbb{S}_{+} \text {and } R \in \mathbb{S}_{++}$ so that $J_{N}$ is convex. We can now formulate the optimal control problem at sampling time $k$ to determine affine control laws with $N$-step prediction horizon
\vspace{-2mm}
\begin{equation}
\begin{array}{cll}
\inf _{H_N} & J_{N}\left(x, H_N\right) &\\
\text { s.t. } & x_{t+1}=A x_{t}+B u_{t}+C w_{t},\quad x_{0}=x,& \forall t \in \left[k,k+N \right]\\
&\sup _{\mathbb{Q} \in \mathbb{B}_{\varepsilon}\left(\widehat{\mathbb{P}}_{N_s}\right)} \mathbb{E}^{\mathbb{Q}}[\ell_j(\xi,H_N)] \le U_j,& \forall j \le N_b ,
\end{array}
\tag{10} \label{eq:DR_MPC}
\end{equation}
where $N_b$ is the number of constraints imposed on state at each sampling time.
This problem with distributionally robust constraints appears to be intractable due to the infinite-dimensional optimization over probability distributions. However, we present a tractable reformulation in the next section.

\section{A TRACTABLE CONVEX CONE PROGRAM REFORMULATION}
In this section we rewrite the distributionally robust control problem \eqref{eq:DR_MPC} with type-$1$ ($p=1$) Wasserstein metric into a finite-dimensional convex cone program leveraging results from robust and convex optimization. After proposing a tractable reformulation, we introduce a practical data-driven Algorithm to handle the disturbed constrained control problem. The control law is computed via cone programming, which enjoy finite sample guarantee, i.e. constraint \eqref{eq:DR_upperbound} is satisfied with regard to a specified level of confidence by collecting a finite number of data points.

We first make some assumptions on the random vector $\xi$ and disturbance $w_k$.
\begin{assumption}[i.i.d. Disturbance]
\label{asp:iid}
We assume that in the discrete-time LTI system \eqref{eq:DLTI}, the disturbance $w_t$ is an i.i.d. random process with  covariance matrix $\Sigma_{w_k}$ and mean $\mu_{w_k}$ for all $t \in \mathbb{N}_0$, which can be computed from data.
\end{assumption}

The i.i.d. random process is a common assumption made in control literature, e.g. \cite{arcari2020dual, coppens2021data}. It assumes a priori that only the first two moments of the random process are acquired as partial distributional information, which can either be estimated or determined a priori \cite{wan2014estimating}.

\begin{assumption}[Moment Assumption \cite{fournier2015rate}]\label{asp:moment}
There exists a positive $\alpha$ such that $
\int_{\Xi} \exp \left(\|\xi\|^{\alpha}\right) \mathbb{Q}(\mathrm{d} \xi)<\infty.
$
\end{assumption}
This assumption trivially holds for a bounded uncertainty set $\Xi$ and finite measure $\mathbb{P}$.

\begin{assumption}[Polytope Uncertainty Set \cite{esfahani2018data}]
\!The space $\mathcal{M}(\Xi)$ of all probability distributions $\mathbb{Q}$ is supported on a polytope $\Xi:=\left\{\xi \in \mathbb{R}^{n_{\xi}}: C_{\xi} \xi \leq d_{\xi}\right\}$.
\end{assumption}

This assumption defines the shape of the uncertainty set. This is a common assumption in the context of robust optimization \cite{LORENZEN2019461} requiring the disturbance not to be infinitely large. This allows the transformation of distribution into linear inequalities.
We then illustrate the equivalent tractable reformulation of the distributionally robust control problem \eqref{eq:disc_quadratic_cost}.
\begin{theorem}[Tractable convex program]
The optimal control problem \eqref{eq:disc_quadratic_cost} with a discounted quadratic cost, distributionally robust constraints within a Wasserstein ball $\mathbb{B}_{\varepsilon}\left(\widehat{\mathbb{P}}_{N_s}\right)$ centered at the empirical distribution $\widehat{P}_{N_s}$ with $N_s$ samples and radius $\varepsilon$ can be reformulated as a cone program \eqref{eq:tractable_MPC} using the equivalent affine control laws from lemma \ref{thm:1} and under Assumptions 1-3.
\begin{equation}
\small \begin{array}{ll} 
\inf _{H_N, \lambda, s_{i}, \gamma_{ijt}} &  \operatorname{Tr}\left\{\left[(\tilde{C}_{x}+B_{x}H_N)^{\top} J_{x}  (\tilde{C}_{x}+B_{x}H_N) + H_N^{\top} J_u H_N \right]\Sigma_w\right\}  + \mu_w^{\top} \left[(\tilde{C}_{x}+B_{x}H_N)^{\top} J_{x} (\tilde{C}_{x}+B_{x}H_N) + H_N^{\top} J_u H_N \right]\mu_w \\
    &\begin{array}{ll}
    & \text { s.t. }   \lambda_j \varepsilon+\frac{1}{N} \sum_{i=1}^{N} s_{ij} \le U_{j} \\
    & b_{tj}+\left\langle a_{tj}, \widehat{\xi}_{i}\right\rangle+\left\langle\gamma_{ijt}, d_{\xi}-C_{\xi} \widehat{\xi}_{i}\right\rangle \leq s_{ij} \\ 
    & \left\|C_{\xi}^{\top} \gamma_{ijt}-a_{tj}\right\|_{*} \leq \lambda_j, \quad \gamma_{ijt} \geq 0\\
    & \quad \forall i \leq N_{s}, \forall j \leq N_b,\forall t \leq N,
    \end{array}
\end{array}
\tag{11} \label{eq:tractable_MPC}
\end{equation}
where $J_{x}:=\operatorname{diag}\left(\operatorname{diag}\left(\beta^{0}, \ldots, \beta^{N-1}\right) \otimes Q, \beta^{N} Q_{f}\right)$, $J_{u}:=\operatorname{diag}\left(\beta^{0}, \ldots, \beta^{N-1}\right) \otimes R$, $a_{tj} = \left[(B_{x} H_N + \tilde{C}_{x} )\right]_{tn_x+m,2:Nn_w+1}$  and $b_{tj} = \left[(B_{x} H_N + \tilde{C}_{x} )\widehat{\xi}_i\right]_{tn_x+m,1}$, where $m \in \mathbb{I}_m$ with the constraint index set $\mathbb{I}_m$ stands for the $m$-th element of the state, on which a constraint is imposed.  $N_s$, $N_b$ and $N$ denote sample number, number of constraints on state at each sampling instance and length of prediction horizon, respectively. $U_j$ is the bound on state. $\widehat{\xi}_{i}$ indicates a data point in the training data set, comprising the disturbance sequence consisted of $N$ sampling time.
\end{theorem}
%
\begin{pf} 
We shall prove the equivalence of the objective function and constraints in \eqref{eq:DR_MPC} and \eqref{eq:tractable_MPC} respectively. 
Applying \eqref{eq:disc_quadratic_cost} with states stacked over the prediction horizon $N$ from \eqref{eq:stacked_state} shows that the objective is a minimax expectation of quadratic cost given a disturbance sequence
\begin{equation}
\begin{split}
\inf _{H_N} \sup _\mathbb{P}\mathbb{E}_{\mathbb{P} \in \mathbb{B}_{\varepsilon}\left(\widehat{\mathbb{P}}_{N_s}\right)}\left\{\tilde{w}_{[k:k+N-1]}^{\top}\left[(\tilde{C}_{x}+B_{x}H_N)^{\top} J_{x} (\tilde{C}_{x} \right.  \left. +B_{x}H_N) + H_N^{\top} J_u H_N \right]\tilde{w}_{[k:k+N-1]}\right\}.
\end{split}
\tag{12} \label{eq:minimax_quadratic_cost}
\end{equation}
Then, under Assumption \ref{asp:iid}, the mean $\mu_{w_{[k:k+N-1]}}$ and covariance matrix $\Sigma_{w_{[k:k+N-1]}}$ of the i.i.d. disturbance sequence are known/computed from the data, the expectation of the quadratic cost is equivalent to the objective function in \eqref{eq:tractable_MPC} according to \cite[THEOREM 1.5]{seber2012linear}.

Given constraints in \eqref{eq:DR_MPC}, representing the worst-case expectation, the linear combination of states is bounded, and we can therefore prove that they are equivalent to constraints in \eqref{eq:tractable_MPC}. This will be proved, without loss of generality, by showing that the upper bound imposed on the $m$-th component of state is satisfied at each sampling time
\begin{equation}
\max( e_m^{\top}\mathbf{e}_t^{\top}x_{[k: k+N]}) \le U_{j}, \forall t \in [1,N], \forall j \le N_b 
\tag{13}
\end{equation}
where $U_{j}$ is a prescribed upper bound on $m$-th component of state and $e_m^{\top}$ selects the $m$-th element of the state.

Given the stacked state represented by the initial state $x_0$ and disturbance sequence in \eqref{eq:stacked_state}, any component of the stacked state within the $N$-step prediction horizon can be written as
\begin{equation}
[x_{[k:k+N]}]_{t  n_x + m} =e_m^{\top}\mathbf{e}_t^{\top} B_{x} H_N(\tilde{C}_y + \tilde{E}_y)\begin{bmatrix} 1\\\xi \end{bmatrix}+\tilde{C}_x \begin{bmatrix} 1\\\xi \end{bmatrix}
\tag{14} \label{eq:stack_state_explicit}
\end{equation}
Thus, we then define a pointwise maximum function consisting of all the $m$-th state elements at each sampling time $t$ for the $j$-th constraint function
\begin{equation}
\ell_j(\xi)=\max_{t \leq N}( [x_{[k:k+N]}]_{t  n_x + m}) = \max _{t \leq N} \left\langle a_{tj}, \xi\right\rangle+b_{tj}
\tag{15}
\end{equation}
where $a_{tj} = \left[(B_{x} H_N + \tilde{C}_{x} )\right]_{tn_x+m,2:Nn_w+1}$  and $b_{tj} = \left[(B_{x} H_N + \tilde{C}_{x} )\xi\right]_{tn_x+m,1}$
to shift the maximum value of $m$-th state entry at each sampling time.

Leveraging the result from \cite[Corollary 5.1]{esfahani2018data}, the distributionally robust constraints in \eqref{eq:DR_MPC} are rewritten into "best-case" constraints 
\begin{equation}
    \inf \lambda_j \varepsilon+\frac{1}{N} \sum_{i=1}^{N} s_{ij} \le U_{j}
    \tag{16}
\end{equation}
along with several additional inequalities. Hence, any feasible solutions of the tractable reformulation \eqref{eq:tractable_MPC} guarantee constraints satisfaction of \eqref{eq:DR_upperbound}. We thus prove the equivalence of the distributionally robust optimization problem \eqref{eq:DR_MPC} and cone program in the form of \eqref{eq:tractable_MPC}.
\end{pf}

\begin{remark}
The lower bound is constructed by setting $\ell_j(\xi)=\max(- [x_{[k:k+N]}]_{t  n_x + m})$ and $U_j$ as negative of the lower bound.
\end{remark}

\begin{remark}
Note that in Section \ref{section:problem_fomulation} we require the state constraints to be polyhedral $C_p x_{[k:k+N]}  \leq D_p$, i.e. only linear combinations of separate state entries.  Therefore, $\ell$ can be chosen to be an affine function of the states, the distributionally robust constraints and \eqref{eq:DR_upperbound}, and can therefore be a function $\ell$ which is affine in the disturbances and initial state. Distributionally robust polyhedral state constraints can then be reformulated into the intersection of linear inequalities as in \eqref{eq:tractable_MPC} by considering joint state constraints of several separate state entries effected by a disturbance sequence.
\end{remark}

We further prove that the  control laws determined by \eqref{eq:tractable_MPC} are able to guarantee joint constraint satisfaction with a finite number of samples. 

\begin{theorem}[finite sample guarantee \cite{esfahani2018data} ]
\label{thm:finite_sample}
\! If Assumption~\ref{asp:moment} (finite moments) holds, and given $H_N$ as the worst-case control law determined via \eqref{eq:tractable_MPC} with ambiguity set $\mathbb{B}_{\varepsilon(N_s,\beta)}\left(\widehat{\mathbb{P}}_{N_s}\right)$ and training data set $\widehat{\Xi}_{N_s}$. Then, for any $p \neq n_wN/2$ the following finite sample guarantee hold with confidence level $1-\beta$ 
\begin{equation}
\mathbb{P}^{N_s}\left\{\mathbb{E}^{\mathbb{P}}[\ell(\xi,H_N)] \leq U \right\} \geq 1-\beta
\tag{17} \label{eq:finite_sample}
\end{equation}
where $\beta \in (0,1)$.
\end{theorem}
\begin{pf}
The finite sample guarantee is the simple consequence of \cite[Theorem 2]{fournier2015rate}. Under Assumption \ref{asp:moment}, the probability that the Wasserstein ball radius does not contain the true probability distribution $\mathbb{P}$ is upper bounded by 
\begin{equation}
\begin{split}
\mathbb{P}^{N_s}\left\{d_{\mathrm{W}}\left(\mathbb{Q}, \widehat{\mathbb{P}}_{N_s}\right) \geq \varepsilon\right\} & \le C \exp{(-cN_s \varepsilon^\kappa)}\mathbb{I}_{\varepsilon\le 1}+ C \exp{(-cN_s \varepsilon^{\alpha/p})}\mathbb{I}_{\varepsilon > 1},
\end{split}
\tag{18} \label{eq:dist_upper_bound}
\end{equation}
where $\kappa(p,\varepsilon) = 2 $ if $p>n_w  N / 2$, and $\kappa(p,\varepsilon) = {d/p}$ if $p \in (0, n_w  N / 2)$. The positive constants $C$ and $c$ depend only on $p,N_s, N,n_w,\alpha$. 

Let $p=1$ for a type-1 Wasserstein metric, we equate the right-hand side of \eqref{eq:dist_upper_bound} to $\beta$ and thus acquire 
\begin{equation}
\small
\varepsilon(N_s, \beta)=\left\{
\begin{array}{cc}\left(\log \left(C \beta^{-1}\right)/(c N_s)\right)^{1 / \kappa } & \text { if } N_s \geq \log \left(C \beta^{-1}\right)/c \\ \left(\log \left(C \beta^{-1}\right)/(c N_s)\right)^{1 / a } & \text { if } N_s < \log \left(C \beta^{-1}\right)/c.\end{array}
\right.
\tag{19} \label{eq:wasserstein_ball_radius_N_beta}
\end{equation}
This directly result in 
\begin{equation}
\mathbb{P}^{N_s}\left\{\mathbb{P} \in \mathbb{B}_{\varepsilon(N_s,\beta)}\left(\widehat{\mathbb{P}}_{N_s}\right) \geq 1-\beta\right\},
\tag{20} \label{eq:prob_chance_constraint}
\end{equation}
where $N_s$ is a finite value. Therefore, $E^{\mathbb{P}}\left[\ell\left(\xi, H_{N}\right)\right] \le \sup _{\mathbb{Q} \in \mathbb{B}_{\varepsilon(N_s,\beta)}\left(\widehat{\mathbb{P}}_{N_s}\right)} \mathbb{E}\left[\ell\left(\xi, H_{N}\right)\right] \leq U$ with probability $1-\beta$.
\end{pf}
\begin{remark}
The proof of Theorem \ref{thm:finite_sample} demonstrates that for any given $\beta$, we can guarantee that  the true distribution is contained within the Wasserstein ball with confidence level $1-\beta$ if we can either collect sufficient samples or expand the ball radius to be large enough.
\end{remark}
\begin{remark}
For convenience of the proof, we required $p \neq n_w N/2$. However, a similar inequality with $\kappa = \varepsilon / \log (2+1 / \varepsilon))^{2}$ holds for $p = n_w  N / 2$.
\end{remark}

After solving problem \eqref{eq:tractable_MPC} we acquire a control law $H_N$ to govern inputs in the succeeding $N$ steps. By following the POB affine control law \eqref{eq:POB_affine}, each input within the prediction horizon can be determined by current and prior purified outputs within the horizon. We can therefore find a policy and then recursively re-solve the conic optimization problem \eqref{eq:tractable_MPC} to update our control law for the subsequent $N_u$ steps. Algorithm \ref{alg:1} illustrates this procedure.
 \vspace{-2mm}
 \begin{algorithm}[H]
 \caption{Distributionally robust MPC}
 \begin{algorithmic}[1]
  \label{alg:1}
   \STATE \textbf{Input: $A$, $B$, $C$, $D$, $E$, $\mu_w$, $\Sigma_w$, $J_x$, $J_u$, $U_j$, $\varepsilon$, $C_\xi$, $d_\xi$ }
    \STATE \textbf{Initialize} $x$, $\tilde{C}_x$, $B_x$, $\mathcal{D}$, $N_s$, $N$, $N_b$, $N_u$, $a_{tj}$, $b_{tj}$, $k=0$
    \REPEAT[every $\Delta T$]
  \IF {($k$ mod $N_u == 0$)}
  \STATE Acquire estimated state $x_k$ and update $\tilde{C}_x$, $a_{tj}$, $b_{tj}$
  \STATE Select samples from $\mathcal{D}$ to formulate $\widehat{\Xi}_{N_s}$
  \STATE Solve the conic optimization problem \eqref{eq:tractable_MPC} for $H_N$
  \STATE Denote the sampling time $t=k$ for control laws update
  \ENDIF
    \STATE Acquire purified observation $v_k$
  \STATE Calculate $u_{k}=h_{k-t}+\sum_{\tau=t}^{k} H_{k-t, \tau-t} v_{\tau}$ by using purified observations collected since time instant $t$ of recent control law update 
  \STATE Store $v_k$, estimate $w_k$ from state evolution and store $w_k$ in
  $\mathcal{D}$
  \STATE update $N_s$ 
  \STATE k = k+1
    \UNTIL END
 \end{algorithmic} 
 \end{algorithm}
 \vspace{-2mm}
\begin{remark}
$N_u$ is the period to update the control law, it can be selected from $\{1,\dots, N-1\}$ depending on the severity of the disturbances and computational cost. To improve the closed-loop performance, a smaller $N_u$ may be selected, while a larger $N_u$ is chosen to reduce computational cost.
\end{remark}
\begin{remark}
$\mathcal{D}$ is the data set containing individual disturbances collected either prior to or during the process. Each sample $\hat{\xi}_i$ fed to problem \eqref{eq:tractable_MPC} is a disturbance sequence of $N$ individual disturbances, i.e. a window of past instances.
It therefore requires that $\mathcal{D}$ contains at least $N$ disturbance data points at initialization. If the disturbance data points are collected in a time sequence, it is possible to replace the oldest data point in the sample $\hat{\xi}_i$ with the newly collected disturbance data point and incorporate one more sample $\hat{\xi}$ into problem \eqref{eq:tractable_MPC} after every $N$ steps. An upper bound may be defined to limit the maximum number of samples incorporated to determine the control laws and cap the computational cost.
\end{remark}

\section{Recursive feasibility}

\begin{assumption}\label{assumption_recur}
We assume that the rank of the matrix $[B_{x}]_{n_x+1:(N+1)n_x, 1:Nn_u}$ is higher than the dimension of the stacked states from $t+1$ to $t+N$ at an arbitrary sampling time $t$, i.e. $\operatorname{rank}([B_{x}]_{n_x+1:(N+1)n_x, 1:Nn_u}) \ge n_x \times N$.
\end{assumption}
\begin{theorem}[Recursive Feasibility]
Let Assumption \ref{assumption_recur} hold. Consider a system given by eq. \eqref{eq:DLTI}, a feasible state trajectory until $t=k$ with $\hat{x}_{[k+1:k+N]} \in \mathbb{D}(x_k)$, and a control law defined by  \eqref{eq:tractable_MPC}. Then, the set of constraints on the states for $t+1$, ${e}_m^{\top} x_i \le U_{j}$ for $i \in [t+1,t+N]$, $\forall j \le N_b$, $\forall m \in \mathbb{I}_m$, are guaranteed for every realization $w_k \in \mathbb{W}$, i.e. the set of feasible states is nonempty $\mathbb{D}(x_{k+1}) \ne \emptyset$.


%
\label{thm:RF}
\end{theorem}

\begin{pf}
We will show that if the optimization problem is feasible at $t = k$ then the bounds on the state elements (i.e. eq. \eqref{eq:DR_upperbound}) hold for the optimization problem \eqref{eq:DR_MPC} at the sampling time $t = k+1$. Without loss of generality, this guarantees recursive feasibility \cite{LOFBERG2012550}.

To prove this, we first reformulate the inequalities in \eqref{eq:DR_upperbound} into an equivalent form in terms of $h_{[0:N-1]}$ and the system matrices, and then prove that there will always exist some $h_{[0:N-1]}$ that is a feasible predicted state trajectory $\hat{x}_{[k+1:k+N]}$.

Let ${e}_m^{\top} x_i \le U_{j} $ be the $j$-th linear constraint on the corresponding $m$-th element of all predicted states within the horizon $N$, where $i \in [t+1,t+N]$ and $j \le N_b$.
Additionally, let $\ell_j(\xi,H_N) = \max_{i \in [t+1:t+N]} \ell_{j,i} = \max_{i \in [t+1:t+N]} {e}_m^{\top} x_i$ be the pointwise maximum of the affine function for the $j$-th constraint.
We can then explicitly represent $\ell_{j,i}$ in terms of $\xi$ and $H_N$.
By reformulating the stacked state \eqref{eq:stacked_state} explicitly in terms of the disturbances, 
we acquire 
\begin{equation}
\ell_{j,i} = {e}_m^{\top} \mathbf{e}_i^{\top}\left(B_{x} H_N(\tilde{C}_y + \tilde{E}_y)\begin{bmatrix} 1\\\xi \end{bmatrix}+\tilde{C}_x \begin{bmatrix} 1\\\xi \end{bmatrix}
\right), \tag{21} \label{eq:pointwise_max_explicit}
\end{equation}
where $\xi = w_{[t:t+N-1]}$, and $\mathbf{e^{\top}_i}$ selects the $i$-th state from the stacked state column vector. Then, by inserting \eqref{eq:pointwise_max_explicit} into the equivalent form  of \eqref{eq:DR_upperbound} \cite{esfahani2018data}, i.e.
\begin{align}
\small
\inf _{\lambda_j \geq 0} \lambda_j \varepsilon+\frac{1}{N} \sum_{i=1}^{N} \sup _{\xi \in \Xi}\left(\ell_j(\xi,H_N)-\lambda\left\|\xi-\widehat{\xi}_{i}\right\|\right) \le U_j \nonumber, \quad \forall j\le N_b 
\tag{22} \label{eq:DR_constraints_reform},
\end{align}
we can reformulate the distributionally robust constraints \eqref{eq:DR_upperbound} into
\begin{align}
\small
\inf _{\lambda_j \geq 0}\sup _{\xi \in \Xi} \min_{\|z_{ijk}\|_{*}\le 1} \lambda_j \varepsilon+\frac{1}{N} \sum_{i=1}^{N} \left(e^{\top}_m x_i- \langle \lambda z_{ijk}, \xi-\widehat{\xi}_{i}\rangle \right) \le U_j, \nonumber \quad \forall i \in [t+1,t+N], \forall j \le N_b, k \le K
\tag{23} \label{eq:bi-dual-reform2}.
\end{align}
Given that there exists a control policy at the sampling time $t = k$ such that the resulting predicted state trajectory $\hat{x}_{[k+1:k+N]}$ satisfies \eqref{eq:bi-dual-reform2} for some $\lambda_j$, $\xi$, $z_{ijk}$. Next, we show the feasibility of optimization problem \eqref{eq:tractable_MPC} at the following sampling time $t = k+1$ by proving that the feasible trajectory $\hat{x}_{[k+1:k+N]}$ at $t = k$ is also feasible at $t = k +1$.
The state $x_i$ is given by $\mathbf{e}_i^{\top}\left(B_{x} H_N(\tilde{C}_y + \tilde{E}_y)\begin{bmatrix} 1\\\xi \end{bmatrix}+\tilde{C}_x \begin{bmatrix} 1\\\xi \end{bmatrix}
\right)$. We let $x_0 = 0$ and $\xi = 0$ to represent the worst case of controllability, i.e. $x_i$ merely dependent on $h_{[0:N-1]}$. The inequality \eqref{eq:bi-dual-reform2} can then be reformulated into 
\begin{equation}
\small
\inf _{\lambda_j \geq 0}\sup _{\xi \in \Xi} \min_{\|z_{ijk}\|_{*}\le 1} \lambda_j \varepsilon+\frac{1}{N} \sum_{i=1}^{N} \left(e^{\top}_j \mathbf{e}_{i}^{T} B_x h_{[0:N-1]}-  \langle \lambda z_{ijk}, \xi-\widehat{\xi}_{i}\rangle \right) \le U_j,
\tag{24} \label{eq:DR_RF_final_form}
\end{equation}
$ \forall i \in [t+1,t+N], \forall j \le N_b, k \le K$.
If the assumption \ref{assumption_recur} holds such that $\hat{x}_{[k+1:k+N]} \in \operatorname{span}([B_{x}]_{n_x+1:(N+1)n_x, 1:Nn_u})$, we can always acquire $\hat{x}_{[k+1:k+N]} = [B_{x}]_{n_x+1:(N+1)n_x, 1:Nn_u} h_{[0:N-1]}$ for some $h_{[0:N-1]}$. Thus, we guarantee that the resulting predicted state trajectory at $t = k+1$ can recover the feasible trajectory $\hat{x}_{[k+1:k+N]}$ at $t=k$ by some control action within the policy. 

\end{pf}
\begin{remark}
If we only impose constraints on $n_s$ separate state elements, then, Assumption \ref{assumption_recur} can be relaxed to the corespondent row of $A^0B,\ldots,A^{N-1}B$ being non-zero and $\operatorname{dim}(h_i) \ge n_s$. In this sense, the predicted state trajectory of a given initial state can be steered to the preceding feasible state trajectory by requiring the constrained state element to be controlled individually at each step. For example, if constraints are imposed on the first state element, we only require the first rows of $A^0B,\ldots,A^{N-1}B$ to be nonzero and one input at each step suffices. This result coincides with the idea of one-step constraints in \cite{7733074}.
\end{remark}

\begin{remark}
We only give a sufficient condition to guarantee the recursive feasibility. Since in the most situations, $H$ would not be eliminated by the worst-case controllability condition. Therefore, we may acquire recursive feasible control policies for a system that does not satisfy the assumption \ref{assumption_recur}. 
\end{remark}
\section{SIMULATIONS AND RESULTS}
In this section, we test the proposed Algorithm \ref{alg:1} on a disturbed mass spring system. We investigate the impact of data samples on the constraint satisfaction, as well as the finite sample guarantee. We also test this framework on an inverted pendulum system and analyse the influence of various Wasserstein ball radii on state trajectories. We outline the simulation specifications and then compare the performance under a number of settings by analysing their constraint and objective values \footnote{The code used for our numerical experiments can be found in \url{https://github.com/OptiMaL-PSE-Lab/Data-driven-distributionally-robust-MPC}.}.

\subsection{Configuration}
Both experiments discretize continuous time models and disturb the systems with sampling period $\Delta T = 0.1$\,s. The prediction horizon is set to $N=5$ and each entry of $w_k$ complies with the random process $3\sin(X)$, where $X \sim \mathcal{N}\left(0, 1\right)$. Therefore, we acquire $C_{\xi} = \operatorname{diag}(1,\dots,-1,\dots)$ and $d_{\xi} = [3,\dots,3] \in \mathbb{R}^{2Nn_w}$. 
The target of both control problems is to steer the disturbed system to the origin whilst satisfying state constraints.

\subsection{Mass spring system}
We consider a mass spring system from \cite{chen1984linear} to illustrate the effectiveness of the proposed Algorithm \ref{alg:1}. 
Maintaining the configuration above, an experiment is conducted with $\varepsilon = 1$ and $N_u = 1$, i.e. the control law is updated at each sampling time.
The weighting matrices are given by $Q = \operatorname{diag}(10,1) $, $Q_f=\operatorname{diag}(15,1)$, $R=1$.
The mass position is directly influenced by 1e-3 times the disturbance and the measurement of the position is noised by another i.i.d. disturbance scaled by 1e-3, i.e. $\small C = \left[\begin{array}{cc}\text{1e-3} & 0 \\ 0 & 0\end{array}\right]$ and $\small E = \left[\begin{array}{cc} 0 & \text{1e-3} \end{array}\right]$. The velocity of the mass is upper bounded by $0.4$\,m/s.

\begin{figure}[thpb]
  \centering
  \includegraphics[width=0.7\textwidth]{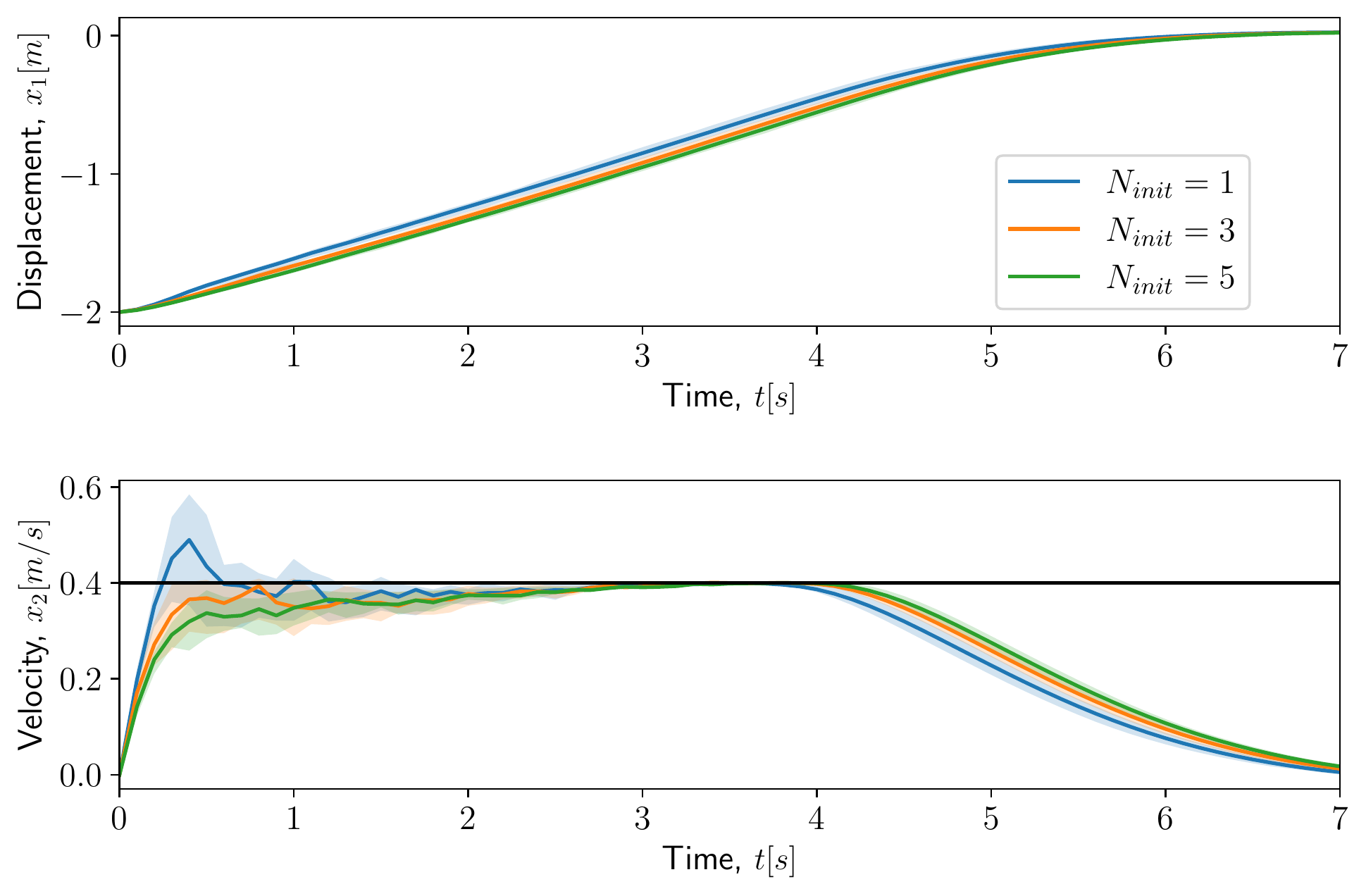}
  \caption{Simulation results of Algorithm \ref{alg:1} from $50$ realizations on a mass spring system. The shaded area denotes the $25$-th to $75$-th trajectory distribution.}
  \label{fig1}
\end{figure}

\begin{figure}[thpb]
  \centering
  \includegraphics[width=0.7\textwidth]{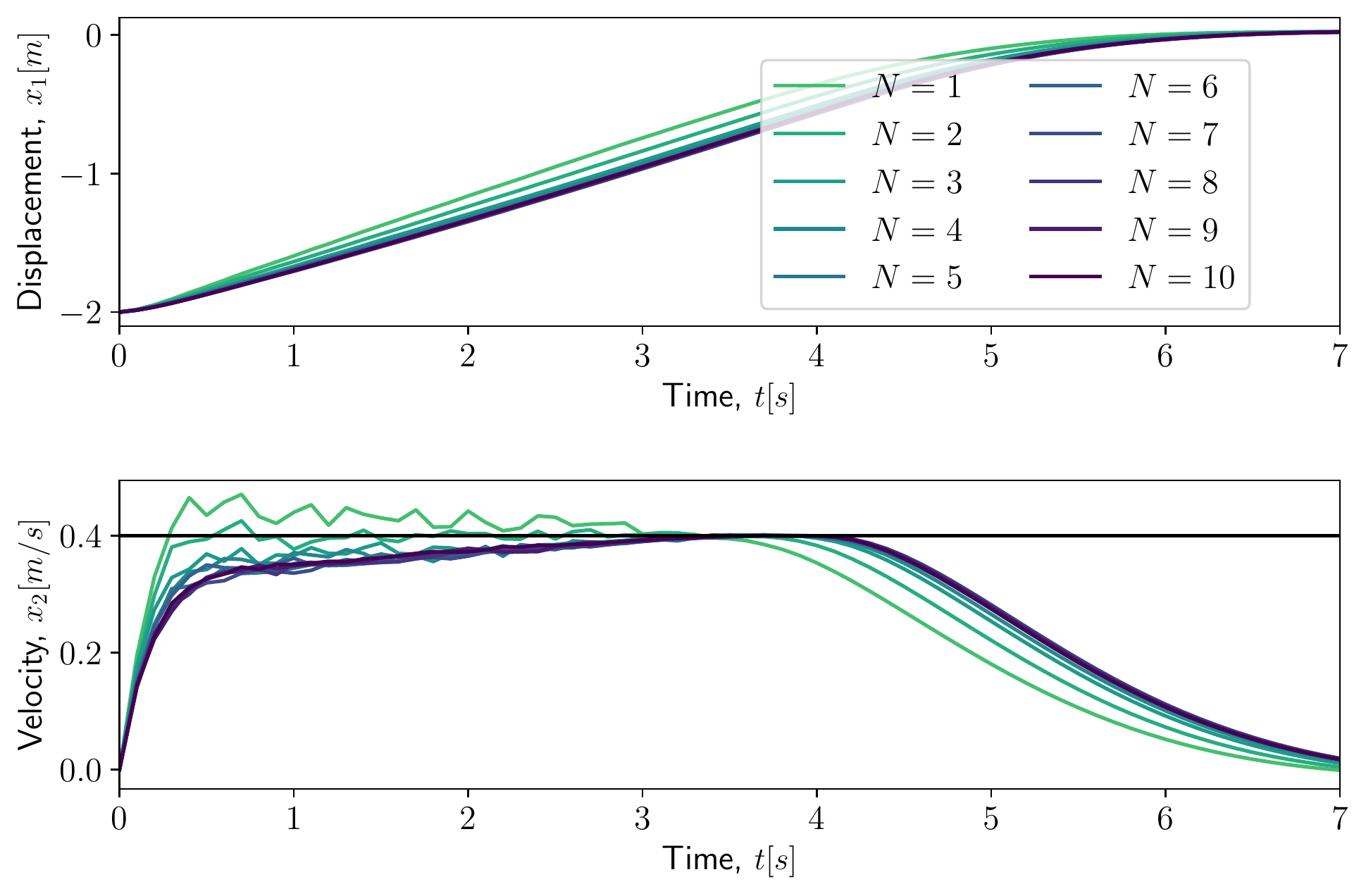}
  \caption{Simulation results of \eqref{eq:tractable_MPC} averaged from $50$ realizations with sample number ranging from $1$ to $10$ on a mass spring system.}
  \label{fig2}
\end{figure}

\begin{figure}[h!]
  \begin{subfigure}[t]{.45\textwidth}
  \centering
    \includegraphics[width=\textwidth]{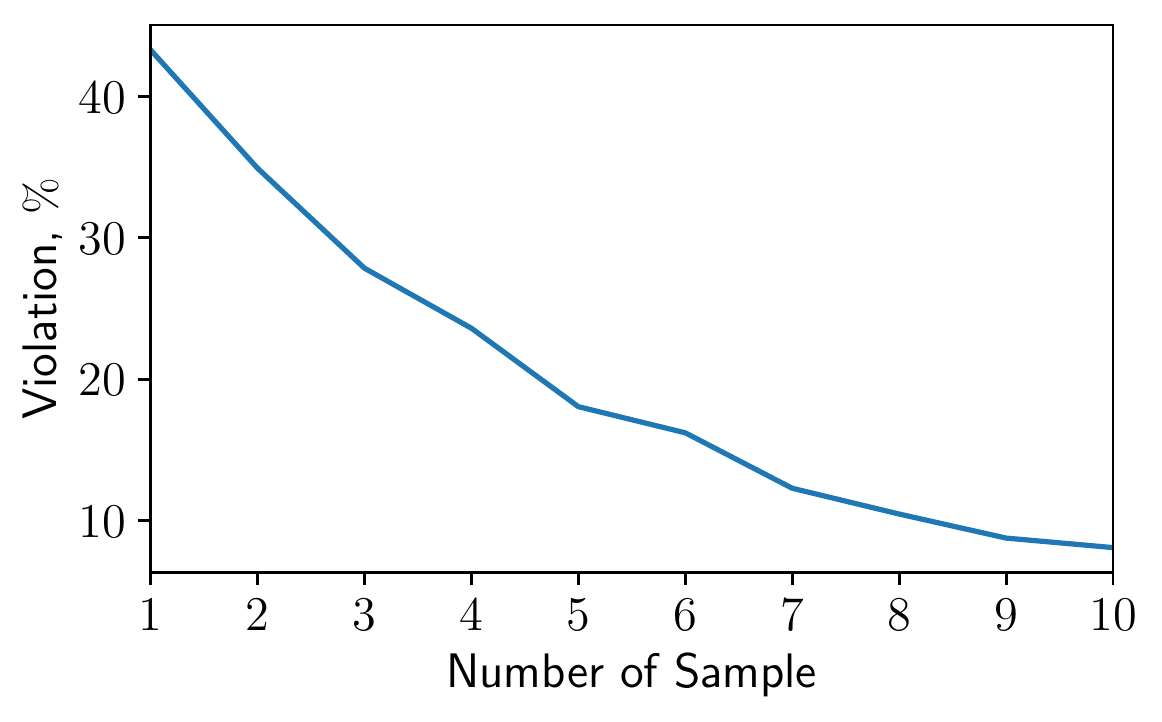}
    \caption{Relation between sample number and constraint violations within first four seconds, averaged from $50$ realizations on the mass spring system.}
      \label{fig3}
  \end{subfigure}\hfill
  \begin{subfigure}[t]{.45\textwidth}
  \centering
    \includegraphics[width=\textwidth]{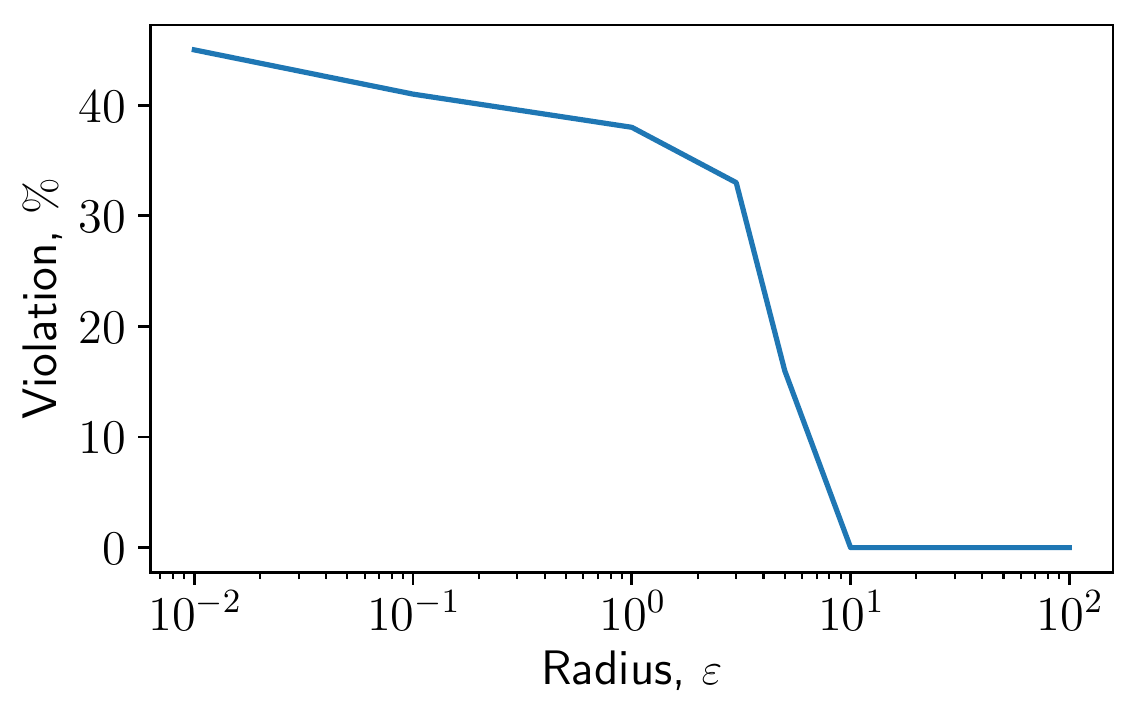}
    \caption{Relation between Wasserstein ball radius and constraint violations within first two seconds, averaged from $10$ realizations on the inverted pendulum system.}
          \label{fig4}
  \end{subfigure}
  \caption{Relation between parameters and constraint violations for the mass spring system and inverted pendulum system. }
\end{figure}

The simulations for Algorithm \ref{alg:1} realize state trajectories with $1,3,5$ collected samples respectively prior to initialization and uses at most 10 recent samples from the data set to avoid large computational costs. As a result, as seen in fig. \ref{fig1}, if we only consider one sample initially, our approach steers the states aggressively. As more samples are used, either collected throughout time or when these samples are provided from the start, the controller steers the state with higher confidence and constraints are satisfied. To analyse our results 50 realizations of the problem are solved.  The shaded area of $25$-th to $75$-th percentiles of simulation initialized with one sample reduces drastically after $1$ second - i.e. after collecting two more samples. Additionally, the $75$-th percentile of trajectory with $N_{init} = 5$ remains feasible during entire simulation, which also manifest higher confidence of constraints satisfaction after acquiring more samples. This is in line with our theoretical results, the number of samples notably increases confidence on constraints satisfaction, and even with few samples the results do not seem over-conservative, we will see this is not the case for enlarging the radius of the Wasserstein ball. Furthermore, the approach is tractable as only a convex cone program is solved. 
In fig. \ref{fig2} we further explore the relation between number of sample and constraint satisfaction. For this example the ball radius is fixed to $1$, the same value as for fig. \ref{fig1}. We simulate the state trajectory with sample numbers ranging from $1$ to $10$, each with $50$ realizations. The control laws are determined at each sampling time with different samples collected prior to initialization. We can read from the figure that the averaged trajectory of $50$ realizations with only $1$ sample tends to violate constraints from the beginning and to oscillate as time increases. In contrast, with large number of samples, constraints are satisfied, this seems to happen for trajectories with sample numbers larger than $5$. Furthermore, results from fig. \ref{fig3} illustrate that the confidence of constraints satisfaction is monotonically increasing along the sample number. 

\subsection{Inverted pendulum}
\begin{figure}[thpb]
  \centering
  \includegraphics[width=0.7\textwidth]{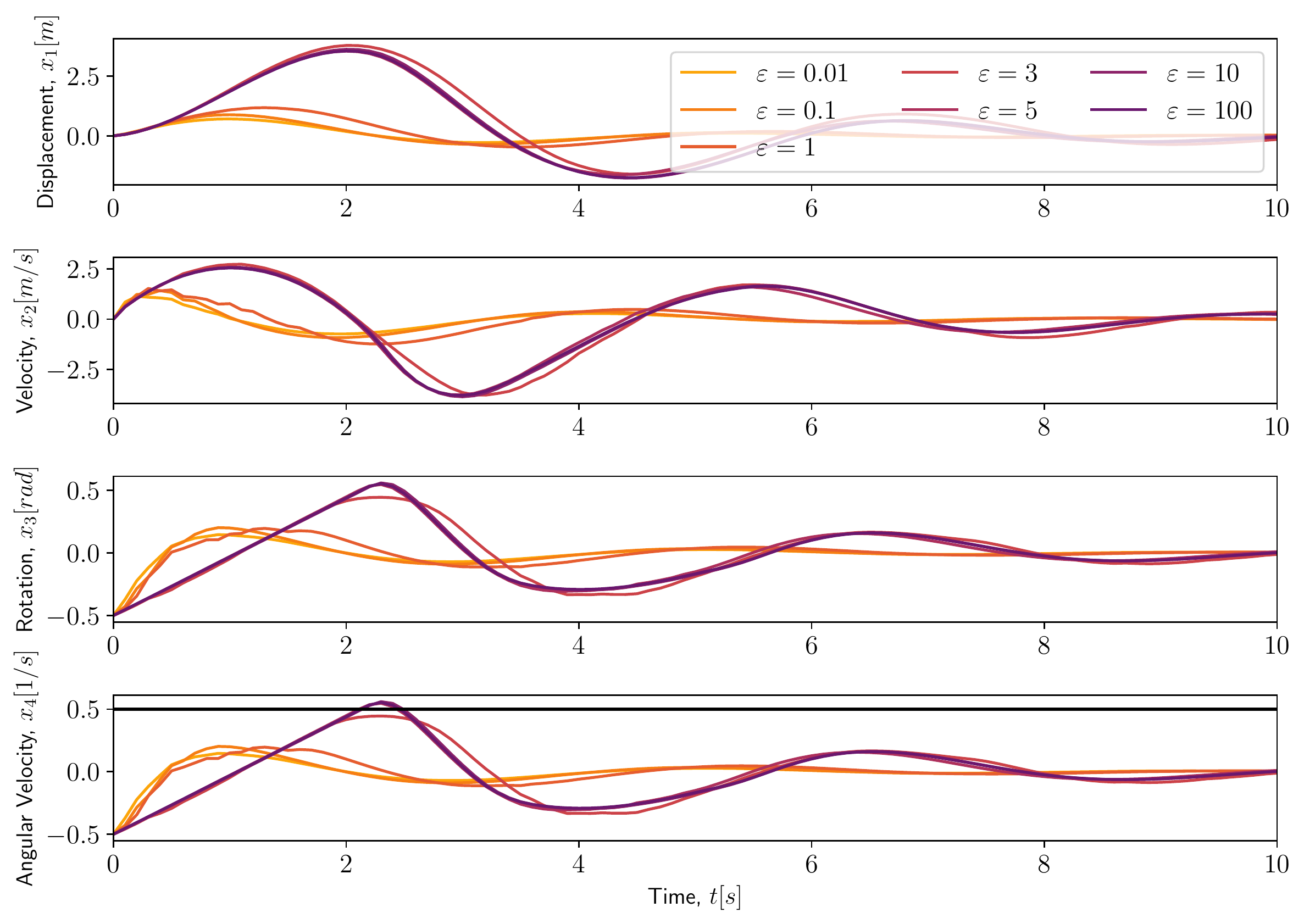}
  \caption{ Simulation results of \eqref{eq:tractable_MPC} averaged from 10 realizations with Wasserstein ball radius ranging from $0.01$ to $100$ on the inverted pendulum.}
  \label{fig5}
\end{figure}

In this case study we illustrate how it is possible to mitigate the constraint violations with a limited number of samples by increasing the Wasserstein ball radius, however, it may lead to conservativeness. We consider an inverted pendulum system represented as the state-space model \cite{singh2017modeling}. Since our interest in this section is to demonstrate the impact of the ball radius, simulations, each with $10$ realizations, are carried out for various ball radii, ranging from $0.01$ to $100$. Control laws are updated at each sampling time with one sample. The  weighting  matrices Q and Qf are $\operatorname{diag}(1000,1,1500,1)$ and $R=1$. The pendulum rotational velocity is disturbed by 1e-2 times the disturbance, whereas the measurement of the pendulum's angular displacement is noised by another i.i.d. disturbance scaled by 1e-2, 
i.e. $\small C^{\top}=\left[\begin{array}{cccc} 0 & 0 & 0 & 0\\ 0 & 0 & 0 & \text{1e-2}  \end{array}\right]$
and $ \small E=\left[\begin{array}{ll}0 &  \text{1e-2}\end{array}\right]$.
The angular velocity is upper bounded by $0.5$\,1/s.
Similar to the first experiment, we use samples collected prior to the initialization to solve the conic problem \eqref{eq:tractable_MPC} at each sampling time in various settings.
As displayed in fig. \ref{fig5}, when the radius is smaller than $1$, the state trajectories violate constraints extensively since the center of the ball is roughly located at one sample's position and very likely this ball does not contain the true distribution. 
If a radius of $3$ or $5$ is used, less constraint violations occur, empirically, constraints are satisfied with probability $65\%$ and $80\%$ respectively. However, if the radius is unnecessarily large, the state trajectories tend to be conservative, e.g. angular velocity is around $4.7$, which is $6\%$ lower than the upper bound. We can see from fig. \ref{fig4} that the confidence of constraint satisfaction increases as the Wasserstein ball expands when the number of samples is fixed.
%
%
%

\vspace{-7.5mm}
\section{CONCLUSIONS}
In this paper, we propose a novel data-driven distributionally robust MPC method for linear systems using the Wasserstein ball. Our approach relies on building an ambiguity set defined by the Wasserstein metric which allows to characterize the uncertainty even when limited information on the probability distributions is available.
In this approach we reformulate the distributionally robust optimal control problem into a tractable convex cone program with finite sample guarantee and propose a practical Algorithm.
Numerical case studies on two systems are conducted to illustrate the effectiveness of the Algorithm and verify the assumptions and theoretical results. Future work includes extending the current approach to incorporate nonlinear dynamics.








\bibliographystyle{cas-model2-names}
\bibliography{ref}

\begin{thebibliography}{38}
\expandafter\ifx\csname natexlab\endcsname\relax\def\natexlab#1{#1}\fi
\providecommand{\url}[1]{\texttt{#1}}
\providecommand{\href}[2]{#2}
\providecommand{\path}[1]{#1}
\providecommand{\DOIprefix}{doi:}
\providecommand{\ArXivprefix}{arXiv:}
\providecommand{\URLprefix}{URL: }
\providecommand{\Pubmedprefix}{pmid:}
\providecommand{\doi}[1]{\href{http://dx.doi.org/#1}{\path{#1}}}
\providecommand{\Pubmed}[1]{\href{pmid:#1}{\path{#1}}}
\providecommand{\bibinfo}[2]{#2}
\ifx\xfnm\relax \def\xfnm[#1]{\unskip,\space#1}\fi
\bibitem[{Arcari et~al.(2020)Arcari, Hewing, Schlichting and
  Zeilinger}]{arcari2020dual}
\bibinfo{author}{Arcari, E.}, \bibinfo{author}{Hewing, L.},
  \bibinfo{author}{Schlichting, M.}, \bibinfo{author}{Zeilinger, M.},
  \bibinfo{year}{2020}.
\newblock \bibinfo{title}{Dual stochastic mpc for systems with parametric and
  structural uncertainty}, in: \bibinfo{booktitle}{Learning for Dynamics and
  Control}, \bibinfo{organization}{PMLR}. pp. \bibinfo{pages}{894--903}.
\bibitem[{Ben-Tal et~al.(2006)Ben-Tal, Boyd and Nemirovski}]{Ben-Tal2006}
\bibinfo{author}{Ben-Tal, A.}, \bibinfo{author}{Boyd, S.},
  \bibinfo{author}{Nemirovski, A.}, \bibinfo{year}{2006}.
\newblock \bibinfo{title}{{Extending scope of robust optimization:
  Comprehensive robust counterparts of uncertain problems}}.
\newblock \bibinfo{journal}{Mathematical Programming} \bibinfo{volume}{107},
  \bibinfo{pages}{63--89}.
\newblock \DOIprefix\doi{10.1007/s10107-005-0679-z}.
\bibitem[{Ben-Tal et~al.(2009)Ben-Tal, El~Ghaoui and
  Nemirovski}]{ben2009robust}
\bibinfo{author}{Ben-Tal, A.}, \bibinfo{author}{El~Ghaoui, L.},
  \bibinfo{author}{Nemirovski, A.}, \bibinfo{year}{2009}.
\newblock \bibinfo{title}{Robust optimization}.
\newblock \bibinfo{publisher}{Princeton university press}.
\bibitem[{Campo and Morari(1987)}]{campo1987robust}
\bibinfo{author}{Campo, P.J.}, \bibinfo{author}{Morari, M.},
  \bibinfo{year}{1987}.
\newblock \bibinfo{title}{Robust model predictive control}, in:
  \bibinfo{booktitle}{1987 American control conference},
  \bibinfo{organization}{IEEE}. pp. \bibinfo{pages}{1021--1026}.
\bibitem[{Chen(1984)}]{chen1984linear}
\bibinfo{author}{Chen, C.T.}, \bibinfo{year}{1984}.
\newblock \bibinfo{title}{Linear system theory and design}. volume
  \bibinfo{volume}{301}.
\newblock \bibinfo{publisher}{Holt, Rinehart and Winston New York}.
\bibitem[{Chen et~al.(2018)Chen, Kuhn and Wiesemann}]{chen2018data}
\bibinfo{author}{Chen, Z.}, \bibinfo{author}{Kuhn, D.},
  \bibinfo{author}{Wiesemann, W.}, \bibinfo{year}{2018}.
\newblock \bibinfo{title}{Data-driven chance constrained programs over
  wasserstein balls}.
\newblock \bibinfo{journal}{arXiv:1809.00210} .
\bibitem[{Coppens and Patrinos(2021)}]{coppens2021data}
\bibinfo{author}{Coppens, P.}, \bibinfo{author}{Patrinos, P.},
  \bibinfo{year}{2021}.
\newblock \bibinfo{title}{Data-driven distributionally robust mpc for
  constrained stochastic systems}.
\newblock \bibinfo{journal}{arXiv preprint arXiv:2103.03006} .
\bibitem[{Dobos et~al.(2009)Dobos, J{\"a}schke, Abonyi and
  Skogestad}]{dobos2009dynamic}
\bibinfo{author}{Dobos, L.}, \bibinfo{author}{J{\"a}schke, J.},
  \bibinfo{author}{Abonyi, J.}, \bibinfo{author}{Skogestad, S.},
  \bibinfo{year}{2009}.
\newblock \bibinfo{title}{Dynamic model and control of heat exchanger networks
  for district heating}.
\newblock \bibinfo{journal}{Hungarian Journal of Industry and Chemistry}
  \bibinfo{volume}{37}.
\bibitem[{Esfahani and Kuhn(2018)}]{esfahani2018data}
\bibinfo{author}{Esfahani, P.M.}, \bibinfo{author}{Kuhn, D.},
  \bibinfo{year}{2018}.
\newblock \bibinfo{title}{Data-driven distributionally robust optimization
  using the wasserstein metric: Performance guarantees and tractable
  reformulations}.
\newblock \bibinfo{journal}{Mathematical Programming} \bibinfo{volume}{171},
  \bibinfo{pages}{115--166}.
\bibitem[{Fournier and Guillin(2015)}]{fournier2015rate}
\bibinfo{author}{Fournier, N.}, \bibinfo{author}{Guillin, A.},
  \bibinfo{year}{2015}.
\newblock \bibinfo{title}{On the rate of convergence in wasserstein distance of
  the empirical measure}.
\newblock \bibinfo{journal}{Probability Theory and Related Fields}
  \bibinfo{volume}{162}, \bibinfo{pages}{707--738}.
\bibitem[{Gao and Kleywegt(2016)}]{gao2016distributionally}
\bibinfo{author}{Gao, R.}, \bibinfo{author}{Kleywegt, A.J.},
  \bibinfo{year}{2016}.
\newblock \bibinfo{title}{Distributionally robust stochastic optimization with
  wasserstein distance}.
\newblock \bibinfo{journal}{arXiv preprint arXiv:1604.02199} .
\bibitem[{Garatti and Campi(2013)}]{garatti2013modulating}
\bibinfo{author}{Garatti, S.}, \bibinfo{author}{Campi, M.C.},
  \bibinfo{year}{2013}.
\newblock \bibinfo{title}{Modulating robustness in control design: Principles
  and algorithms}.
\newblock \bibinfo{journal}{IEEE Control Systems Magazine}
  \bibinfo{volume}{33}, \bibinfo{pages}{36--51}.
\bibitem[{Guo et~al.(2018)Guo, Baker, Dall’Anese, Hu and
  Summers}]{guo2018data}
\bibinfo{author}{Guo, Y.}, \bibinfo{author}{Baker, K.},
  \bibinfo{author}{Dall’Anese, E.}, \bibinfo{author}{Hu, Z.},
  \bibinfo{author}{Summers, T.H.}, \bibinfo{year}{2018}.
\newblock \bibinfo{title}{Data-based distributionally robust stochastic optimal
  power flow—part i: Methodologies}.
\newblock \bibinfo{journal}{IEEE Transactions on Power Systems}
  \bibinfo{volume}{34}, \bibinfo{pages}{1483--1492}.
\bibitem[{Hewing et~al.(2020)Hewing, Wabersich and
  Zeilinger}]{HEWING2020109095}
\bibinfo{author}{Hewing, L.}, \bibinfo{author}{Wabersich, K.P.},
  \bibinfo{author}{Zeilinger, M.N.}, \bibinfo{year}{2020}.
\newblock \bibinfo{title}{Recursively feasible stochastic model predictive
  control using indirect feedback}.
\newblock \bibinfo{journal}{Automatica} \bibinfo{volume}{119},
  \bibinfo{pages}{109095}.
\newblock \DOIprefix\doi{https://doi.org/10.1016/j.automatica.2020.109095}.
\bibitem[{Houska and M.E(2019)}]{RMPC_book_Houska}
\bibinfo{author}{Houska, B.}, \bibinfo{author}{M.E, V.}, \bibinfo{year}{2019}.
\newblock \bibinfo{title}{Robust optimization for MPC}.
\bibitem[{Lorenzen et~al.(2019)Lorenzen, Cannon and
  Allgöwer}]{LORENZEN2019461}
\bibinfo{author}{Lorenzen, M.}, \bibinfo{author}{Cannon, M.},
  \bibinfo{author}{Allgöwer, F.}, \bibinfo{year}{2019}.
\newblock \bibinfo{title}{Robust mpc with recursive model update}.
\newblock \bibinfo{journal}{Automatica} \bibinfo{volume}{103},
  \bibinfo{pages}{461--471}.
\newblock \DOIprefix\doi{https://doi.org/10.1016/j.automatica.2019.02.023}.
\bibitem[{Lorenzen et~al.(2017)Lorenzen, Dabbene, Tempo and
  Allgöwer}]{7733074}
\bibinfo{author}{Lorenzen, M.}, \bibinfo{author}{Dabbene, F.},
  \bibinfo{author}{Tempo, R.}, \bibinfo{author}{Allgöwer, F.},
  \bibinfo{year}{2017}.
\newblock \bibinfo{title}{Constraint-tightening and stability in stochastic
  model predictive control}.
\newblock \bibinfo{journal}{IEEE Transactions on Automatic Control}
  \bibinfo{volume}{62}, \bibinfo{pages}{3165--3177}.
\newblock \DOIprefix\doi{10.1109/TAC.2016.2625048}.
\bibitem[{Lu et~al.(2020)Lu, Lee and You}]{You2020}
\bibinfo{author}{Lu, S.}, \bibinfo{author}{Lee, J.H.}, \bibinfo{author}{You,
  F.}, \bibinfo{year}{2020}.
\newblock \bibinfo{title}{Soft-constrained model predictive control based on
  data-driven distributionally robust optimization}.
\newblock \bibinfo{journal}{AIChE Journal} \bibinfo{volume}{66},
  \bibinfo{pages}{e16546}.
\bibitem[{Lucia et~al.(2012)Lucia, Finkler, Basak and Engell}]{LUCIA201269}
\bibinfo{author}{Lucia, S.}, \bibinfo{author}{Finkler, T.},
  \bibinfo{author}{Basak, D.}, \bibinfo{author}{Engell, S.},
  \bibinfo{year}{2012}.
\newblock \bibinfo{title}{A new robust nmpc scheme and its application to a
  semi-batch reactor example*}.
\newblock \bibinfo{journal}{IFAC Proceedings Volumes} \bibinfo{volume}{45},
  \bibinfo{pages}{69--74}.
\newblock \DOIprefix\doi{https://doi.org/10.3182/20120710-4-SG-2026.00035}.
  \bibinfo{note}{8th IFAC Symposium on Advanced Control of Chemical Processes}.
\bibitem[{Lucia and Paulen(2014)}]{LUCIA20141904}
\bibinfo{author}{Lucia, S.}, \bibinfo{author}{Paulen, R.},
  \bibinfo{year}{2014}.
\newblock \bibinfo{title}{Robust nonlinear model predictive control with
  reduction of uncertainty via robust optimal experiment design}.
\newblock \bibinfo{journal}{IFAC Proceedings Volumes} \bibinfo{volume}{47},
  \bibinfo{pages}{1904--1909}.
\newblock \DOIprefix\doi{https://doi.org/10.3182/20140824-6-ZA-1003.02332}.
\bibitem[{Luo and Mehrotra(2019)}]{luo2019decomposition}
\bibinfo{author}{Luo, F.}, \bibinfo{author}{Mehrotra, S.},
  \bibinfo{year}{2019}.
\newblock \bibinfo{title}{Decomposition algorithm for distributionally robust
  optimization using wasserstein metric with an application to a class of
  regression models}.
\newblock \bibinfo{journal}{European Journal of Operational Research}
  \bibinfo{volume}{278}, \bibinfo{pages}{20--35}.
\bibitem[{Löfberg(2012)}]{LOFBERG2012550}
\bibinfo{author}{Löfberg, J.}, \bibinfo{year}{2012}.
\newblock \bibinfo{title}{Oops! i cannot do it again: Testing for recursive
  feasibility in mpc}.
\newblock \bibinfo{journal}{Automatica} \bibinfo{volume}{48},
  \bibinfo{pages}{550--555}.
\bibitem[{Mayne et~al.(2011)Mayne, Kerrigan, Van~Wyk and
  Falugi}]{mayne2011tube}
\bibinfo{author}{Mayne, D.Q.}, \bibinfo{author}{Kerrigan, E.C.},
  \bibinfo{author}{Van~Wyk, E.}, \bibinfo{author}{Falugi, P.},
  \bibinfo{year}{2011}.
\newblock \bibinfo{title}{Tube-based robust nonlinear model predictive
  control}.
\newblock \bibinfo{journal}{International Journal of Robust and Nonlinear
  Control} \bibinfo{volume}{21}, \bibinfo{pages}{1341--1353}.
\bibitem[{Mesbah et~al.(2014)Mesbah, Streif, Findeisen and
  Braatz}]{mesbah2014stochastic}
\bibinfo{author}{Mesbah, A.}, \bibinfo{author}{Streif, S.},
  \bibinfo{author}{Findeisen, R.}, \bibinfo{author}{Braatz, R.D.},
  \bibinfo{year}{2014}.
\newblock \bibinfo{title}{Stochastic nonlinear model predictive control with
  probabilistic constraints}, in: \bibinfo{booktitle}{2014 American control
  conference}, \bibinfo{organization}{IEEE}. pp. \bibinfo{pages}{2413--2419}.
\bibitem[{Nilim and El~Ghaoui(2005)}]{nilim2005robust}
\bibinfo{author}{Nilim, A.}, \bibinfo{author}{El~Ghaoui, L.},
  \bibinfo{year}{2005}.
\newblock \bibinfo{title}{Robust control of markov decision processes with
  uncertain transition matrices}.
\newblock \bibinfo{journal}{Operations Research} \bibinfo{volume}{53},
  \bibinfo{pages}{780--798}.
\bibitem[{{Petsagkourakis} et~al.(2020){Petsagkourakis}, {Heath}, {Carrasco}
  and {Theodoropoulos}}]{9145693}
\bibinfo{author}{{Petsagkourakis}, P.}, \bibinfo{author}{{Heath}, W.P.},
  \bibinfo{author}{{Carrasco}, J.}, \bibinfo{author}{{Theodoropoulos}, C.},
  \bibinfo{year}{2020}.
\newblock \bibinfo{title}{Robust stability of barrier-based model predictive
  control}.
\newblock \bibinfo{journal}{IEEE Transactions on Automatic Control} ,
  \bibinfo{pages}{1--1}\DOIprefix\doi{10.1109/TAC.2020.3010770}.
\bibitem[{Piccoli and Rossi(2014)}]{piccoli2014generalized}
\bibinfo{author}{Piccoli, B.}, \bibinfo{author}{Rossi, F.},
  \bibinfo{year}{2014}.
\newblock \bibinfo{title}{Generalized wasserstein distance and its application
  to transport equations with source}.
\newblock \bibinfo{journal}{Archive for Rational Mechanics and Analysis}
  \bibinfo{volume}{211}, \bibinfo{pages}{335--358}.
\bibitem[{Rahimian and Mehrotra(2019)}]{rahimian2019distributionally}
\bibinfo{author}{Rahimian, H.}, \bibinfo{author}{Mehrotra, S.},
  \bibinfo{year}{2019}.
\newblock \bibinfo{title}{Distributionally robust optimization: A review}.
\newblock \bibinfo{journal}{arXiv preprint arXiv:1908.05659} .
\bibitem[{Rawlings et~al.(2017)Rawlings, Mayne and Diehl}]{MPC_book}
\bibinfo{author}{Rawlings, J.}, \bibinfo{author}{Mayne, D.},
  \bibinfo{author}{Diehl, M.}, \bibinfo{year}{2017}.
\newblock \bibinfo{title}{Model Predictive Control: Theory, Computation, and
  Design}.
\bibitem[{Schwarm and Nikolaou(1999)}]{Schwarm1999}
\bibinfo{author}{Schwarm, A.T.}, \bibinfo{author}{Nikolaou, M.},
  \bibinfo{year}{1999}.
\newblock \bibinfo{title}{Chance-constrained model predictive control}.
\newblock \bibinfo{journal}{AIChE Journal} \bibinfo{volume}{45},
  \bibinfo{pages}{1743--1752}.
\newblock \DOIprefix\doi{https://doi.org/10.1002/aic.690450811}.
\bibitem[{Seber and Lee(2012)}]{seber2012linear}
\bibinfo{author}{Seber, G.A.}, \bibinfo{author}{Lee, A.J.},
  \bibinfo{year}{2012}.
\newblock \bibinfo{title}{Linear regression analysis}. volume
  \bibinfo{volume}{329}.
\newblock \bibinfo{publisher}{John Wiley \& Sons}.
\bibitem[{Singh and Singla(2017)}]{singh2017modeling}
\bibinfo{author}{Singh, G.}, \bibinfo{author}{Singla, A.},
  \bibinfo{year}{2017}.
\newblock \bibinfo{title}{Modeling, analysis and control of a single stage
  linear inverted pendulum}, in: \bibinfo{booktitle}{2017 IEEE International
  Conference on Power, Control, Signals and Instrumentation Engineering
  (ICPCSI)}, \bibinfo{organization}{IEEE}. pp. \bibinfo{pages}{2728--2733}.
\bibitem[{Touretzky and Baldea(2014)}]{TOURETZKY20141292}
\bibinfo{author}{Touretzky, C.R.}, \bibinfo{author}{Baldea, M.},
  \bibinfo{year}{2014}.
\newblock \bibinfo{title}{Integrating scheduling and control for economic mpc
  of buildings with energy storage}.
\newblock \bibinfo{journal}{Journal of Process Control} \bibinfo{volume}{24},
  \bibinfo{pages}{1292--1300}.
\newblock \DOIprefix\doi{https://doi.org/10.1016/j.jprocont.2014.04.015}.
  \bibinfo{note}{economic nonlinear model predictive control}.
\bibitem[{Van~Parys et~al.(2015)Van~Parys, Kuhn, Goulart and
  Morari}]{van2015distributionally}
\bibinfo{author}{Van~Parys, B.P.}, \bibinfo{author}{Kuhn, D.},
  \bibinfo{author}{Goulart, P.J.}, \bibinfo{author}{Morari, M.},
  \bibinfo{year}{2015}.
\newblock \bibinfo{title}{Distributionally robust control of constrained
  stochastic systems}.
\newblock \bibinfo{journal}{IEEE Transactions on Automatic Control}
  \bibinfo{volume}{61}, \bibinfo{pages}{430--442}.
\bibitem[{Wan et~al.(2014)Wan, Wang, Liu and Tong}]{wan2014estimating}
\bibinfo{author}{Wan, X.}, \bibinfo{author}{Wang, W.}, \bibinfo{author}{Liu,
  J.}, \bibinfo{author}{Tong, T.}, \bibinfo{year}{2014}.
\newblock \bibinfo{title}{Estimating the sample mean and standard deviation
  from the sample size, median, range and/or interquartile range}.
\newblock \bibinfo{journal}{BMC medical research methodology}
  \bibinfo{volume}{14}, \bibinfo{pages}{1--13}.
\bibitem[{Wiesemann et~al.(2014)Wiesemann, Kuhn and
  Sim}]{wiesemann2014distributionally}
\bibinfo{author}{Wiesemann, W.}, \bibinfo{author}{Kuhn, D.},
  \bibinfo{author}{Sim, M.}, \bibinfo{year}{2014}.
\newblock \bibinfo{title}{Distributionally robust convex optimization}.
\newblock \bibinfo{journal}{Operations Research} \bibinfo{volume}{62},
  \bibinfo{pages}{1358--1376}.
\bibitem[{Yang(2019)}]{yang2019data}
\bibinfo{author}{Yang, I.}, \bibinfo{year}{2019}.
\newblock \bibinfo{title}{Data-driven distributionally robust stochastic
  control of energy storage for wind power ramp management using the
  wasserstein metric}.
\newblock \bibinfo{journal}{Energies} \bibinfo{volume}{12},
  \bibinfo{pages}{4577}.
\bibitem[{{Yang}(2020)}]{9222209}
\bibinfo{author}{{Yang}, I.}, \bibinfo{year}{2020}.
\newblock \bibinfo{title}{Wasserstein distributionally robust stochastic
  control: A data-driven approach}.
\newblock \bibinfo{journal}{IEEE Transactions on Automatic Control}
  \DOIprefix\doi{10.1109/TAC.2020.3030884}.

\end{thebibliography}

\end{document}